\documentclass[11pt]{article}

\usepackage{graphicx}

\usepackage{amssymb,subfigure}
\usepackage{amsmath}

\newtheorem{teo}{Theorem}

\newtheorem{df}{Definition}

\renewcommand{\r}{{\mathbb R}}
\newcommand{\z} {{\mathbb Z}}
\newcommand{\cn} {{\mathbb C}}
\newcommand{\n} {{\mathbb N}}
\newcommand{\e} {\varepsilon}
\newcommand{\be}{\begin{equation}}
\newcommand{\ee}{\end{equation}}
\newcommand{\bega}{\begin{gather}}
\newcommand{\enga}{\end{gather}}

\begin{document}

\title{An inequality for a periodic uncertainty constant
}

\author{Elena A. Lebedeva\footnote{Mathematics and Mechanics Faculty, Saint Petersburg State University,
Universitetsky prospekt, 28, Peterhof,  Saint Petersburg,
 198504, Russia; Saint Petersburg State Polytechnical University,
 Polytechnicheskay 29,  Saint Petersburg, 195251, Russia}}

\date{ealebedeva2004@gmail.com}

\maketitle


\begin{abstract}
 An inequality refining the lower bound for a periodic (Breitenberger) uncertainty constant is proved for a wide class of functions.   
A connection of uncertainty constants for periodic and non-periodic functions is extended to this class. A particular minimization problem for a non-periodic (Heisenberg) uncertainty constant is studied.   
\end{abstract}



\section{Introduction}
\label{intr}

The Breitenberger uncertainty constant (UC) is commonly used as a measure of localization for periodic functions.  It was introduced in 1985 by Brei\-tenberger in \cite{B}. It can be derived from a general operator "position-momentum" approach as it is discussed in \cite{FolSit}. The Breitenberger UC has a deep connection with the classical Heisenberg UC, which characterizes localization of functions on the real line. There exists a universal lower bound for both UCs (the uncertainty principle). It equals  $1/2$ (see  chosen normalization in Sec. \ref{note}). It is well known that the least value is attained on the Gaussian function in the real line case and there is no such  function in the periodic case. 
At the same time, in \cite{Bat} Battle proves a number of inequalities specifying the lower bound of the  Heisenberg UC for wavelets. In particular, it is proved that if a wavelet $\psi^0\in L_2(\r)$ has a zero frequency centre 
$c(\widehat{\psi^0}):=\int_{\r}\xi |\widehat{\psi^0}(\xi)|^2 \, d\xi/(\int_{\r}|\widehat{\psi^0}(\xi)|^2 \, d\xi)=0,$ 
then the Heisenberg
UC is greater or equal to $3/2$ (see \cite[Theorem 1.4]{Bat}).  
  
  The main contribution of this paper is an inequality refining the lower bound of the Breitenberger UC for a wide class of  sequences of periodic functions (Theorem \ref{main}). This result is somewhat analogous to Battle's result mentioned above. Given a sequence of periodic functions $\psi_j,$ $j \in \z_+,$ the conditions $|(\psi'_j,\,\psi_j)|\leq C \|\psi_j\|^2$ and $\lim_{j\to \infty} q_j \widehat{\psi}_j(k) / \|\psi_j\| =0$ (see (\ref{cond3}) and (\ref{cond1}) in Theorem \ref{main}) correspond to a zero frequency centre 
$c(\widehat{\psi^0})=0$ and the wavelet admissibility condition $\widehat{\psi^0}(0)=0$ respectively.  The rest of restrictions (\ref{cond2}), (\ref{cond4})--(\ref{cond6}) in Theorem \ref{main}  mean some ``regularity'' of the sequence $\psi_j.$
In \cite{prqurase03},  the following formula connecting  UCs for periodic ($UC_B$) and non-periodic ($UC_H$) functions is obtained 
$
\lim_{j\to\infty} UC_B(\psi^p_j)=UC_H(\psi_0),
$
where  $\psi^p_{j}(x):=2^{j/2}\sum_{n \in \mathbb{Z}} \psi^0(2^j(x+2 \pi n)),$ $j \in\z_+$.
In Step 3 of the proof of Theorem \ref{main},  we generalize this formula and suggest a new proof of this fact. 
In Remark 1 -- Remark 3, we discuss which classes of periodic wavelet sequences satisfy the conditions of Theorem \ref{main}. 
We also study one particular minimization problem for the Heisenberg UC connected with Battle's result mentioned above (Theorem \ref{nominfunc}).  
If the result of Theorem \ref{nominfunc} had been wrong, it would have been possible to give another proof of Theorem \ref{main}.  

While there are sufficiently many results specifying the lower and the upper bounds of the Heisenberg UC \cite{Balan,Bat,FolSit,GB,GL,L08,L07,L12,L11,N1} and the upper bound of the Breitenberger UC \cite{LebPres14,NW,PQ,R05,Se},  
to our knowledge, there are not actually any results concerting to an estimation of the lower bound for the Breitenberger UC in the literature.
  
  This work also has the following motivation. In \cite{LebPres14},  a family of periodic Parseval wavelet frames is constructed. The family has optimal time-frequency localization (the Breitenberger UC tends to $1/2$) 
with respect to a family parameter, and it has the best currently known localization (the Breitenberger UC tends to $3/2$) with respect to a multiresolution analysis parameter. 
In \cite{LebPres14}, the conjecture was formulated: the Breitenberger UC is greater than $3/2$ for any periodic wavelet sequence $(\psi_j)_{j\in\z_+}$ such that $(\psi_j',\,\psi_j)_{L_{2,2\pi}}=0$. Theorem \ref{main} of this paper proves the conjecture for a wide class of  sequences of periodic functions under a milder restriction $(\psi_j',\,\psi_j)_{L_{2,2\pi}} \leq C \|\psi_j\|^2_{L_{2,2\pi}}$. So the family constructed in \cite{LebPres14} has optimal localization with respect to both parameters within the class of functions considered in Theorem \ref{main}.

\section{Notations and auxiliary results}
\label{note}
Let $L_{2,2 \pi}$ be the space of all $2\pi$-periodic square-integrable complex-valued functions, with
inner product $(\cdot,\cdot)$ given by 
$
(f,\,g):= (2\pi)^{-1} \int_{-\pi}^{\pi} f(x)\overline{g(x)}\,\mathrm{d}x
$
for any $f,g \in L_{2,2\pi},$ and norm $\|\cdot\|:=\sqrt{(\cdot,\,\cdot)}.$ 
The Fourier series of a function 
$
f \in L_{2,2\pi}
$
is defined by
$\sum_{k \in \mathbb{Z}}\widehat{f}(k) \mathrm{e}^{ \mathrm{i} k x},$
where its Fourier coefficient is defined by
$
\widehat{f}(k) = (2\pi)^{-1} \int_{-\pi}^{\pi} f(x)\mathrm{e}^{- \mathrm{i} k x}\,\mathrm{d}x.
$

Let $L_{2}(\r)$ be the space of all square-integrable complex-valued functions, with
inner product $(\cdot,\cdot)$ given by 
$
(f,\,g):= (2\pi)^{-1} \int_{\r} f(x)\overline{g(x)}\,\mathrm{d}x
$
for any $f,g \in L_2(\r),$ and norm $\|\cdot\|:=\sqrt{(\cdot,\,\cdot)}.$ 
The Fourier transform of a function 
$
f \in L_{2}(\r)
$
is defined by
$
\widehat{f}(\xi):= (2\pi)^{-1} \int_{\r}f(x) \mathrm{e}^{- \mathrm{i} \xi x}\,\mathrm{d}x.
$

Let us recall the definitions of the  UCs and the uncertainty principles.
 
 \begin{df}[\cite{H}]
 \label{Huc}
\texttt{The Heisenberg UC} of $f \in L_2(\mathbb{R})$ is the functional 
$UC_H(f):=\Delta(f)\Delta(\widehat{f})$ such that
$$
\Delta^2(f):=
\|f\|^{-2}\|(\cdot-c(f))f\|^{2}, \ \ \ \ \ \ \ 
c(f):=
\|f\|^{-2}(\cdot f,\,f),
$$
where 
$\Delta(f),$ $\Delta(\widehat{f}),$ $c(f),$ and  $c(\widehat{f})$
are called \texttt{time variance, frequency variance, time centre,} and \texttt{frequency centre} respectively.
\end{df}

It is clear that the time variance can be rewritten as 
$$
\Delta^2(f)=\frac{\|\cdot f\|^2}{\|f\|^2}-\frac{(\cdot f,\,f)^2}{\|f\|^4}
$$ 
Using  elementary properties of the Fourier transform,  we rewrite the frequency variance as 
$$
\Delta^2(\widehat{f})=\frac{\| i f'\|^2}{\|f\|^2} - \frac{(i f',\,f)^2}{\|f\|^4}
$$ 
(See \cite{prqurase03} Lemmas 1 and 2, where this trick is explained in detail).

\begin{teo}[\cite{H, FolSit}; the Heisenberg uncertainty principle]
Let $f \in L_2(\mathbb{R})$, then $UC_H(f)\geq 1/2,$ and the equality is attained iff $f$ is the Gaussian function. 
\end{teo}

\begin{teo} [\cite{Bat}, p.137; refinement of the Heisenberg uncertainty principle] 
\label{Battle}
If $f \in L_2(\r),$ $\ c(\widehat{f})$ $=0$, and $\int_{\r} f = 0,$  
then $UC_H(f)\geq 3/2.$
\end{teo}

\begin{df}[\cite{B}]
\label{uc}
Let $f =\sum_{k \in \mathbb{Z}} c_k \mathrm{e}^{ \mathrm{i} k \cdot}\in L_{2,2\pi}.$  
\texttt{ The first trigonometric moment} is defined as
$$ 
\tau(f):=\frac{1}{2 \pi} \int_{-\pi}^{\pi} \mathrm{e}^{ \mathrm{i} x} |f(x)|^2\, \mathrm{d}x =
 \sum_{k \in \mathbb{Z}} c_{k-1} \overline{c_{k}}.
$$
\texttt{The angular variance} of the function $f$ is defined by  
$$
{\rm var_A }(f):= \frac{\left(\sum_{k \in \mathbb{Z}}|c_k|^2\right)^2}{
\left|\sum_{k \in \mathbb{Z}}c_{k-1} \bar{c}_{k}\right|^2}-1
=
\frac{\|f\|^4}{|\tau(f)|^2}-1.
$$
\texttt{The frequency variance} of the function $f$ is defined by  
$$
{\rm var_F }(f):= \frac{\sum_{k \in \mathbb{Z}}k^2 |c_k|^2}{\sum_{k \in \mathbb{Z}}|c_k|^2}-
\frac{\left(\sum_{k \in \mathbb{Z}}k|c_k|^2\right)^2}{\left(\sum_{k \in \mathbb{Z}}|c_k|^2\right)^2}
=
\frac{\|f'\|^2}{\|f\|^2}-\frac{(i f',\, f)^2}{\|f\|^4}.
$$
The quantity 
$
UC_B(\{c_k\}):=UC_B(f):=\sqrt{\mathrm{var_A}(f)\mathrm{var_F}(f)}
$
is called \texttt{the  Breitenberger (periodic) UC}.
\end{df}

We consider also two additional terms to characterize the first trigonometric moment (In another form, they are introduced in \cite[Lemma 3]{prqurase03}). Namely,  by definition, put
\be
\label{defAB}
A(f):=\frac12\sum_{k\in \z} |c_{k-1}-c_k|^2,\ \  \ \ 
B(f):=\frac12 \sum_{k\in \z}(c_{k-1}-c_k)(\overline{c_{k-1}}+\overline{c_{k}})
\ee
 for $f(x)=\sum_{k \in \mathbb{Z}} c_k \mathrm{e}^{ \mathrm{i} k x}\in L_{2,2\pi}.$
 It is clear that 
 \be
 \label{AB}
 A(f) = \sum_{k\in \z} |c_k|^2 - \Re \left(\sum_{k\in \z} c_{k-1}\overline{c_k}\right) = \|f\|^2 - \Re (\tau(f)), 
 \ \ \ \ 
 B(f) = i \Im \left(\sum_{k\in \z} c_{k-1}\overline{c_k}\right) = i \Im (\tau(f)). 
 \ee

\begin{teo}[\cite{B, PQ}; the Breitenberger  uncertainty principle]
\label{UC}
Let $f \in L_{2,2\pi}$,   $f(x)\neq C 
\mathrm{e}^{ \mathrm{i} k x},$ $C \in \mathbb{R},$ $k \in \mathbb{Z}$. Then
 $UC_B(f) > 1/2$  and there is no function such that  $UC_B(f) = 1/2.$
\end{teo}


Now, we recall the notion of a tight frame.
Let $H$ be a separable Hilbert space. If there exists a constant $A>0$
such that for any $f \in H$ the following equality holds
$
\sum_{n=1}^{\infty} \left|(f,\,f_n)\right|^2 = A \|f\|^2,
$ 
 then the sequence $(f_n)_{n \in \mathbb{N}}$ is called \texttt{a tight frame} for $H.$ In the case $A=1$, a tight frame is called \texttt{ a Parseval frame}.
In addition, if $\|f_n\|=1$ for all $n \in \mathbb{N}$,  then a Parseval frame  forms an orthonormal basis.

We are especially interested to get a refinement of the Breitenberger uncertainty principle for the case of periodic wavelet sequences. For our purposes, it is sufficient  to consider wavelet systems with one wavelet generator. We recall the basic notions.   
In the sequel,  we use the following notation 
$
f_{j,k}(x):=f_j(x- 2 \pi 2^{-j} k)
$
for a function $f_j \in L_{2,2\pi}.$
Consider  functions   $\varphi_0,\,\psi_j \in L_{2,2\pi},$ $j\in\z_+$. If the set 
 $\Psi:=\left\{\varphi_0,  \psi_{j,k} : \ j\in\z_+,\ k=0,\dots,2^j-1 \right\}$
forms a tight frame (or a basis) for $L_{2,2\pi}$ then $\Psi$ is said to be \texttt{a periodic tight wavelet frame 
(or a periodic wavelet basis)} for $L_{2,2\pi}.$

\begin{teo}[\cite{GT1}; the unitary extension principle for a periodic setting]
\label{UEP}
Let $\varphi_j \in L_{2,2\pi},$ $j\in\z_+$,
be  a sequence of $2\pi$-periodic functions such that 
\begin{equation}\label{con1}
     \lim_{j \to \infty}2^{j/2} \widehat{\varphi}_j(k) = 1, \ \ \ \ \ \ \  k \in \z.
\end{equation}
Let $\mu^j_k \in \cn,$ $j\in\z_+,$ $k\in\z$, be a two-parameter sequence    such that $\mu^j_{k+2^j}=\mu^j_{k},$
and
\begin{equation}
\label{con2}
	\widehat{\varphi}_{j}(k)=\mu^{j+1}_k \widehat{\varphi}_{j+1}(k). 
\end{equation}
Let $\psi_j,$ $j\in\z_+,$ 
be a sequence of $2\pi$-periodic functions defined using Fourier coefficients
\begin{equation}
\label{con3}
		{\widehat{\psi}}_{j}(k)=\lambda^{j+1}_{k} \widehat{\varphi}_{j+1}(k),
\end{equation}
where $\lambda^j_{k}\in\cn$, $\lambda^j_{k+2^j}=\lambda^j_{k}$,  and 
\begin{equation}
\label{con4}
	\left(
	\begin{array}{cc}
	\mu ^j_k &  \mu ^j_{k+2^{j-1}} \\
	 \lambda^j_{k} & \lambda^j_{k+2^{j-1}}
	\end{array}
	\right)
	\left(
	\begin{array}{cc}
	\overline{\mu} ^j_k &  \overline{\lambda}^j_{k} \\
	\overline{\mu} ^j_{k+2^{j-1}} & \overline{\lambda}^j_{k+2^{j-1}}
	\end{array}
	\right)
	=
	\left(
	\begin{array}{cc}
	2 & 0 \\
	0 & 2
	\end{array}
	\right).
\end{equation}
Then the family
$\Psi : =\left\{\varphi_0,  \psi_{j,k} : \ j \in\z_+,\ k=0,\dots,2^j-1\right\}$
forms a Parseval wavelet frame for $L_{2,2\pi}.$ 
\end{teo}
The sequences $(\varphi_j)_{j\in\z_+},$ $(\psi_j)_{j\in\z_+},$ $(\mu^j_k)_{k\in\z},$ and $(\lambda^j_k)_{k\in\z}$ are 
called \texttt{a scaling sequence, a wavelet sequence, a scaling mask and a wave\-let mask}   respectively.

A periodic wavelet system can be constructed starting with a scaling mask. 
Namely, let $\nu^{j}_{k}$ be  a sequence given by  $\nu^{j}_{k}=\nu^{j}_{k+2^j}$.	We define 
$\widehat{\xi}_{j}(k):=\prod_{r=j+1}^{\infty}\nu^{r}_{k}.$ 
If the above infinite products converge, then the scaling sequence, scaling mask, wavelet mask, and wavelet sequence are defined 
	 respectively as
\begin{gather}
        \widehat{\varphi_{j}}(k):=2^{-j/2}\widehat{\xi}_{j}(k),\qquad\quad 
	\mu^{j}_k:=\sqrt{2} \nu^{j}_k,\ 
	\notag\\
	\lambda^{j}_k:=e^{2\pi i 2^{-j}k}\mu^{j}_{k+2^{j-1}},\qquad \quad
	{\widehat{\psi}}_{j}(k):=\lambda^{j+1}_{k} \widehat{\varphi}_{j+1}(k).
\notag
\end{gather}

\section{Main result}

In the following theorem we prove an inequality for the Breitenberger UC for a wide class of sequences of periodic functions.  

\begin{teo}
\label{main}
Let $\psi_j \in L_{2,2\pi},$ $j\in\n$ be periodic functions such that  
\begin{gather}
\label{cond1}
\lim_{j\to \infty} q_j \widehat{\psi}_j(k) / \|\psi_j\| =0 \mbox{ for } |k| \leq M(C), \\
\label{cond2}
\lim_{j\to\infty} q_j^{-2} A(\psi'_j) /  \|\psi_j\|^2= 0,\\  
\label{cond3}
|(\psi'_j,\,\psi_j)|\leq C \|\psi_j\|^2,\\
\label{cond4}
q_j^{-2}\|\psi'_j\|^2\leq C \|\psi_j\|^2, \\ 
\label{cond5}
q_j^{2} A(\psi_j)\leq C \|\psi_j\|^2, \\ 
\label{cond6}
q_{j} \left|B(\psi_j)\right|\leq C \|\psi_j\|^2
\end{gather}
where $M(C):=2 (C +  C \sqrt{2C}/3 +1/6),$
 $C>0$ is an absolute constant,  
and $q_j \rightarrow +\infty$ as $j \to +\infty.$
If
$\lim_{j \to \infty}UC_B(\psi_j)$ exists,
 then
 \be 
\lim_{j \to \infty}UC_B(\psi_j)\geq 3/2.
\ee
\end{teo}

\textbf{Proof.} 
\textbf{Step 1.}
The Breitenberger UC  is a homogeneous functional of degree zero, 
that is $UC(\alpha f) = UC(f)$ for  $\alpha \in \cn \setminus\{0\}.$ So in the sequel, we consider the functions $\psi_j/\|\psi_j\|$ instead of $\psi_j.$ However, to avoid the fussiness of notations we keep the former name for the function $\psi_j.$  
Thus we consider a sequence $(\psi_j)_{j \in \z_+},$ where $\|\psi_j\|=1$, and conditions (\ref{cond1}) -- (\ref{cond6}) take the form
$$
\lim_{j\to \infty} q_j \widehat{\psi}_j(k)=0 \mbox{ for }  |k| \leq M(C),\ \ \ 
\lim_{j\to\infty} q_j^{-2} A(\psi'_j)  = 0,\ \ \  
|(\psi'_j,\,\psi_j)|\leq C ,
$$
$$
q_j^{-2}\|\psi'_j\|\leq C, \ \ \ 
q_j^{2} A(\psi_j)\leq C,\ \ \ 
q_j B(\psi_j)\leq C.
$$

\textbf{Step 2.}
Let us consider auxiliary functions $\psi_{\ast j}\in L_{2,2\pi}$ such that 
\be
\label{condpsi0}
\widehat{\psi_{\ast j}}(k)=\left\{
\begin{array}{ll}
\widehat{\psi_{j}}(k), & |k| > M(C), \\
0, & |k| \leq M(C).
\end{array}
\right.
\ee
It follows from (\ref{cond1}) that 
$$
\|\psi_j\|^2-\|\psi_{\ast j}\|^2 = \sum_{|k|\leq M(C)}|\widehat{\psi_{j}}(k)|^2 = q_j^{-2} o(1),
\ \ \ \ 
\tau(\psi_j)-\tau(\psi_{\ast j}) = \sum_{|k|\leq M(C)}\widehat{\psi_{j}}(k-1)\overline{\widehat{\psi_{j}}(k)} = q_j^{-2} o(1),
$$
$$
(\psi'_j,\, \psi_j)-(\psi'_{\ast j},\, \psi_{\ast j}) = i  \sum_{|k|\leq M(C)}k|\widehat{\psi_{j}}(k)|^2 = q_j^{-2} o(1),
\ \ \ \ 
\|\psi'_j\|^2-\|\psi'_{\ast j}\|^2 = \sum_{|k|\leq M(C)}k^2|\widehat{\psi_{j}}(k)|^2 = q_j^{-2} o(1)
$$
as $j\to\infty.$
It is straightforward to see that 
\be
\label{UCBpsi0j}
 UC_B(\psi_{\ast j})-UC_B(\psi_{j}) \rightarrow 0 \mbox{ as } j \to \infty.
\ee
Indeed,
it follows from (\ref{AB}), $\|\psi_j\|=1$, and (\ref{cond5}) that 
$$
1\geq |\tau(\psi_j)|^2 = (\|\psi_j\|^2 - A(\psi_j))^2 + ( i B(\psi_j))^2=
(1 - A(\psi_j))^2 + ( i B(\psi_j))^2 \geq 1-2 A(\psi_j) \geq 1- 2 C q_j^{-2}.
$$
So 
$$
0\leq {\rm var_A }(\psi_j) = \frac{1}{|\tau(\psi_j)|^2} - 1 \leq 
\frac{2C q_j^{-2} }{1-2C q_j^{-2} }.
$$
Therefore, $q_j^2 {\rm var_A }(\psi_j)$ is bounded as $j \to \infty.$

By definition of ${\rm var_F }$, and (\ref{cond4}) we obtain the boundedness of $q_j^{-2}{\rm var_F }(\psi_{\ast j})$
$$
0\leq q_j^{-2}{\rm var_F }(\psi_{\ast j}) = 
q_j^{-2}  \left(\frac{\|(\psi_{\ast j})'\|^2}{\|\psi_{\ast j}\|^2}-\frac{(i (\psi_{\ast j})',\, \psi_{\ast j})^2}{\|\psi_{\ast j}\|^4}\right) \leq
q_j^{-2} \frac{\|(\psi_{j})'\|^2+q_j^{-2} o(1)}{1-q_j^{-2} o(1)} \leq 
\frac{C+q_j^{-4} o(1)}{1-q_j^{-2} o(1)}.
$$
We rewrite now $UC_B(\psi_{\ast j})-UC_B(\psi_{j})$ in a standard form
$$
 \left(q_j^{2}{\rm var_A }(\psi_{\ast j})-q_j^{2}{\rm var_A }(\psi_{ j})\right)q_j^{-2}{\rm var_F }(\psi_{\ast j}) + \left(q_j^{-2}{\rm var_F }(\psi_{\ast j})-q_j^{-2}{\rm var_F }(\psi_{ j})\right)q_j^{2}{\rm var_A }(\psi_{ j}),
$$
and we estimate 
$q_j^{2} \left({\rm var_A }(\psi_{\ast j})-{\rm var_A }(\psi_{ j})\right)$ and
$q_j^{-2}\left({\rm var_F }(\psi_{\ast j})-{\rm var_F }(\psi_{ j})\right).$
Using the above estimates for $\|\psi_j\|^2-\|\psi_{\ast j}\|^2$ and 
$\tau(\psi_j)-\tau(\psi_{\ast j})$  we get
$$
q_j^{2} \left|{\rm var_A }(\psi_{\ast j})-{\rm var_A }(\psi_{ j})\right| = 
q_j^{2} \left|\frac{\|\psi_{\ast j}\|^4}{|\tau(\psi_{\ast j})|^2}-\frac{\|\psi_{j}\|^4}{|\tau(\psi_{j})|^2}\right|$$
$$
\leq
q_j^{2} \frac{
\left|\|\psi_{\ast j}\|^2-\|\psi_{ j}\|^2\right|
\left(\|\psi_{\ast j}\|^2+\|\psi_{ j}\|^2\right)
|\tau(\psi_{\ast j})|^2
+
\left||\tau(\psi_{\ast j})|-|\tau(\psi_{j})|\right| 
\left(|\tau(\psi_{\ast j})|+|\tau(\psi_{j})|\right) 
\|\psi_{j}\|^4
}{|\tau(\psi_{\ast j})|^2|\tau(\psi_{ j})|^2}
$$
$$
\leq
q_j^2
\frac{
q_j^{-2}o(1) (1+q_j^{-2}o(1))+2 q_j^{-2}o(1)
}
{(1-2q_j^{-2})(1-2(C+o(1))q_j^{-2})}.
$$
So  $q_j^{2} \left|{\rm var_A }(\psi_{\ast j})-{\rm var_A }(\psi_{ j})\right| \to 0$ as $j \to \infty.$
Analogously, using the above estimates for $\|\psi'_j\|^2-\|\psi'_{\ast j}\|^2$ and 
$(\psi'_j,\, \psi_j)-(\psi'_{\ast j},\, \psi_{\ast j})$  we have
$$
q_j^{-2}\left|{\rm var_F }(\psi_{\ast j})-{\rm var_F }(\psi_{ j})\right|
=
q_j^{-2}\left|
\frac{\|\psi'_{\ast j}\|^2}{\|\psi_{\ast j}\|^2}+\frac{(\psi'_{\ast j},\psi_{\ast j})^2}{\|\psi_{\ast j}\|^4}- \frac{\|\psi'_{ j}\|^2}{\|\psi_{ j}\|^2}-\frac{(\psi'_{j},\psi_{j})^2}{\|\psi_{j}\|^4}
\right|
$$
$$
\leq
q_j^{-2}
\frac{
\left|\|\psi'_{\ast j}\|^2-\|\psi'_{ j}\|^2\right| 
\|\psi_{\ast j}\|^2
+ 
\left|\|\psi_{\ast j}\|^2-\|\psi_{ j}\|^2\right|
\|\psi'_{j}\|^2 }
{\|\psi_{\ast j}\|^2\|\psi_{ j}\|^2}
$$
$$
+
q_j^{-2}
\frac{
\left|(\psi'_{\ast j},\psi_{\ast j})-(\psi'_{ j},\psi_{ j})\right|
\left|(\psi'_{\ast j},\psi_{\ast j})+(\psi'_{ j},\psi_{ j})\right|
\|\psi_{\ast j}\|^4 + 
\left|\|\psi_{\ast j}\|^2-\|\psi_{ j}\|^2\right|
\left(\|\psi_{\ast j}\|^2+\|\psi_{ j}\|^2\right)
\left|(\psi'_{ j},\psi_{ j})^2\right|
}{\|\psi_{\ast j}\|^4\|\psi_{ j}\|^4}
$$
$$
\leq 
q_j^{-2}\frac{q_j^{-2}o(1)(1+q_j^{-2}o(1))+q_j^{-2}o(1) C q_j^{2}}{1-q_j^{-2}o(1)}
$$
$$
+
q_j^{-2}\frac{q_j^{-2}o(1)(2 C+q_j^{-2}o(1))(1-q_j^{-2}o(1)) + q_j^{-2}o(1)(1+q_j^{-2}o(1))C}{1-q_j^{-2}o(1)}.
$$
So $q_j^{-2}\left({\rm var_F }(\psi_{\ast j})-{\rm var_F }(\psi_{ j})\right) \to 0$ as $j\to \infty$, and (\ref{UCBpsi0j}) is proved.

\textbf{Step 3.} 
Let us introduce auxiliary functions $f_{\ast j},$ $j\in\n$ such that 
\begin{gather}
\label{conf1}
q_j^{-1/2}\widehat{f_{\ast j}}(q_j^{-1}k) = \widehat{\psi_{\ast j}}(k), 
\\
\label{conf2} 
 \widehat{f_{\ast j}} \mbox{ is linear on any interval } [q_j^{-1}(k-1),\,q_j^{-1}k],
 \end{gather}
 where $k \in \z.$ Since $\psi_{\ast j} \in L_{2,2\pi}$, it follows that $f_{\ast j} \in L_2(\r)$ for $j\in\z.$ We claim that 
 \be
 \label{UCHUCB}
 UC_B(\psi_{\ast j}) - UC_H(f_{\ast j}) \rightarrow 0 \mbox{ as } j\to \infty.
 \ee    
Indeed, we estimate differences 

\begin{enumerate}
	\item $\|\psi_{\ast j}\|^2_{L_{2,2\pi}} - \|f_{\ast j}\|^2_{L_2(\r)}, \ \ \ $
	\item $q_j^{-1}(\psi'_{\ast j},\,\psi_{\ast j})_{L_{2,2\pi}}-(f'_{\ast j},\,f_{\ast j})_{L_2(\r)}, \ \ \ $
	\item $q_j^{-2}\|\psi'_{\ast j}\|^2_{L_{2,2\pi}} - \|f'_{\ast j}\|^2_{L_2(\r)}, $
	\item $2\cdot q_j^{2} A(\psi_{\ast j}) - \|\cdot f_{\ast j}\|^2_{L_2(\r)}, \ \ \ $
	\item $q_j B(\psi_{\ast j}) - i(\cdot f_{\ast j},\,f_{\ast j})_{L_2(\r)}.$
\end{enumerate}

1.
For the first difference, consequently using the definition of $f_{\ast j}$ (\ref{conf2}), (\ref{AB}), and (\ref{cond5}) we have
$$
\|\psi_{\ast j}\|^2_{L_{2,2\pi}} - \|f_{\ast j}\|^2_{L_2(\r)} = 
\sum_{k\in\z} \bigl|\widehat{\psi_{\ast j}}(k)\bigr|^2 -\int_{\r}\bigl|\widehat{f_{\ast j}}(\xi)\bigr|^2\,d\xi=
\sum_{k\in\z} \bigl|\widehat{\psi_{\ast j}}(k)\bigr|^2 -\sum_{k\in\z} \int_{q_j^{-1}(k-1)}^{q_j^{-1}k}
\bigl|\widehat{f_{\ast j}}(\xi)\bigr|^2\,d\xi
$$
$$
=\sum_{k\in\z} \bigl|\widehat{\psi_{\ast j}}(k)\bigr|^2 -\sum_{k\in\z} \int_{k-1}^{k}  \left|\left(\widehat{\psi_{\ast j}}(k-1)-\widehat{\psi_{\ast j}}(k)\right)(\xi-k)+\widehat{\psi_{\ast j}}(k)\right|^2\,d\xi
=\sum_{k\in\z} \bigl|\widehat{\psi_{\ast j}}(k)\bigr|^2 
$$
$$
-\frac13 \left(\sum_{k\in\z} \bigl|\widehat{\psi_{\ast j}}(k)\bigr|^2 +\sum_{k\in\z} \bigl|\widehat{\psi_{\ast j}}(k-1)\bigr|^2+\Re \left(\sum_{k\in\z} \widehat{\psi_{\ast j}}(k-1)\overline{\widehat{\psi_{\ast j}}(k)}\right)  \right)= \frac13 \left(\|\psi_{\ast j}\|^2_{L_{2,2\pi}}-\Re(\tau(\psi_{\ast j}))\right)
$$
$$
=\frac13 A(\psi_{\ast j})\leq \frac{C}{3} q_j^{-2}
$$
So, 
\be
\label{estnorm}
\|\psi_{\ast j}\|^2_{L_{2,2\pi}} - \|f_{\ast j}\|^2_{L_2(\r)}\leq  \frac{C}{3} q_j^{-2} \to 0 
\mbox{ as } j \to \infty.
\ee
2.
Similarly to the previous calculation, we obtain
$$
\left|q_j^{-1}(\psi'_{\ast j},\,\psi_{\ast j})_{L_{2,2\pi}}-(f'_{\ast j},\,f_{\ast j})_{L_2(\r)}\right|=
\left|q_j^{-1} i \sum_{k\in\z} k \bigl|\widehat{\psi_{\ast j}}(k)\bigr|^2 - i \int_{\r} \xi \bigl|\widehat{f_{\ast j}}(\xi)\bigr|^2\,d\xi\right|
$$
$$
= q_j^{-1}  \left|\sum_{k\in\z} k \bigl|\widehat{\psi_{\ast j}}(k)\bigr|^2 -\sum_{k\in\z} \int_{k-1}^{k} \xi \left|\left(\widehat{\psi_{\ast j}}(k-1)-\widehat{\psi_{\ast j}}(k)\right)(\xi-k)+\widehat{\psi_{\ast j}}(k)\right|^2\,d\xi \right|
$$
$$
\leq \frac{q_j^{-1}}{3}  \left| \sum_{k\in\z} k \left(\bigl|\widehat{\psi_{\ast j}}(k)\bigr|^2-\Re \bigl( \widehat{\psi_{\ast j}}(k-1)\overline{\widehat{\psi_{\ast j}}(k)}\bigr)\right)\right|
+\frac{q_j^{-1}}{12}\left|\Re  \left(\sum_{k\in\z} \widehat{\psi_{\ast j}}(k-1)\overline{\widehat{\psi_{\ast j}}(k)}\right)
\right|
$$
For the first sum, the Cauchy inequality, (\ref{cond4}), and (\ref{cond5}) yield
$$
\frac{q_j^{-1}}{3}\left| \sum_{k\in\z} k \left(\bigl|\widehat{\psi_{\ast j}}(k)\bigr|^2-\Re \bigl( \widehat{\psi_{\ast j}}(k-1)\overline{\widehat{\psi_{\ast j}}(k)}\bigr)\right)\right|
\leq \frac{q_j^{-1}}{3}
\left(\sum_{k\in\z} k^2 \bigl|\widehat{\psi_{\ast j}}(k)\bigr|^2\right)^{1/2}
$$
$$
\times \left(\sum_{k\in\z}\bigl|\widehat{\psi_{\ast j}}(k)-\widehat{\psi_{\ast j}}(k-1)\bigr|^2\right)^{1/2}
 \leq \frac{q_j^{-1}}{3} \|\psi'_{\ast j}\|_{L_{2,2\pi}}\left(2 A(\psi_{\ast j})\right)^{1/2}\leq \frac{2^{1/2}}{3} C^{3/2} q_j^{-1}. 
$$
For the second one, using (\ref{AB}) and (\ref{cond5}), we have
$$
\frac{q_j^{-1}}{12}\left|\Re  \left(\sum_{k\in\z} \widehat{\psi_{\ast j}}(k-1)\overline{\widehat{\psi_{\ast j}}(k)}\right)
\right|=\frac{q_j^{-1}}{12}\left|\Re \bigl(\tau(\psi_{\ast j})\bigr)\right|=
\frac{q_j^{-1}}{12} \left|\|\psi_{\ast j}\|^2-A(\psi_{\ast j})\right|\leq \frac{q_j^{-1}}{6}.
$$
Therefore,
\be
\label{estcentre}
\left|q_j^{-1}(\psi'_{\ast j},\,\psi_{\ast j})_{L_{2,2\pi}}-(f'_{\ast j},\,f_{\ast j})_{L_2(\r)}\right|
\leq \left(\frac{C\sqrt{2C}}{3} + \frac{1}{6}\right)q_j^{-1} \to 0 
\mbox{ as } j \to \infty.
\ee

3. The next estimations are analogous to the previous ones, so we omit details 
$$
\left|q_j^{-2}\|\psi'_{\ast j}\|^2_{L_{2,2\pi}} - \|f'_{\ast j}\|^2_{L_2(\r)}\right| = 
\left|q_j^{-2}  \sum_{k\in\z} k^2 \bigl|\widehat{\psi_{\ast j}}(k)\bigr|^2 - \int_{\r} \xi^2 \bigl|\widehat{f_{\ast j}}(\xi)\bigr|^2\,d\xi\right|
$$
$$
= q_j^{-2} \left| 
\frac16 \sum_{k\in \z}
 \bigl| k \widehat{\psi_{\ast j}}(k) - (k-1) \widehat{\psi_{\ast j}}(k-1) \bigr|^2 -
 \frac{1}{30} \left(2\sum_{k\in \z}
 \bigl|  \widehat{\psi_{\ast j}}(k) \bigr|^2 + 
3 \Re \left(\sum_{k\in \z}\widehat{\psi_{\ast j}}(k-1)\overline{\widehat{\psi_{\ast j}}(k)} \right) 
  \right)
 \right|
$$
$$
= q_j^{-2} \left| 
\frac16  A(\psi'_{\ast j})- \frac{1}{30} \left(2\|\psi_{\ast j}\|^2   +3 \Re\bigl(\tau(\psi_{\ast j})\bigr) \right)  
\right|=
q_j^{-2} \left| 
\frac16  A(\psi'_{\ast j})- \frac{1}{30} \left(5 \|\psi_{\ast j}\|^2-3 A(\psi_{\ast j}) \right)  
\right|
$$
Then, by (\ref{cond2}), (\ref{cond5}), we conclude that 
$$
q_j^{-2}\|\psi'_{\ast j}\|^2_{L_{2,2\pi}} - \|f'_{\ast j}\|^2_{L_2(\r)} \to 0 \mbox{ as } 
j\to \infty.
$$

Due to piecewise linearity of the functions $\widehat{f_{\ast j}}$,  the differences in item 4 and item 5 are equal to $0$ for all $j \in \n.$ 

4. Indeed, exploiting the definitions of $A$ (\ref{defAB}) and $f_{\ast j}$ (\ref{conf1}), we immediately get
$$
2\cdot q_j^{2} A(\psi_{\ast j}) - \|\cdot f_{\ast j}\|^2_{L_2(\r)} = 
q_j^{2} \sum_{k\in\z} \left|\widehat{\psi_{\ast j}}(k-1) - \widehat{\psi_{\ast j}}(k)\right|^2 - 
\int_{\r} \left|\widehat{f_{\ast j}}'(\xi)\right|^2\,d\xi 
$$
$$
= q_{j} \sum_{k\in\z} \left|\widehat{f_{\ast j}}\bigl(q_j^{-1}(k-1)\bigr) - \widehat{f_{\ast j}}\bigl(q_j^{-1} k\bigr)\right|^2 - 
\int_{\r} \left|\widehat{f_{\ast j}}'(\xi)\right|^2\,d\xi 
$$
$$
= \sum_{k\in \z} \int_{q_j^{-1}(k-1)}^{q_j^{-1} k}\left(\left|\frac{\widehat{f_{\ast j}}\bigl(q_j^{-1}(k-1)\bigr) - \widehat{f_{\ast j}}\bigl(q_j^{-1} k\bigr)}{q_j^{-1}}\right|^2 -\left|\widehat{f_{\ast j}}'(\xi)\right|^2  \right)\,d\xi =0.
$$

5. The definitions of $B$ (\ref{defAB}) and $f_{\ast j}$ (\ref{conf1})  yields 
$$
q_j B(\psi_{\ast j}) - i(\cdot f_{\ast j},\,f_{\ast j})_{L_2(\r)} = 
 q_j B(\psi_{\ast j}) + 2\pi (\widehat{f_{\ast j}}',\,\widehat{f_{\ast j}})_{L_2(\r)} 
 $$
 $$
 =
 q_{j}\sum_{k\in\z} \frac{\overline{\widehat{\psi_{\ast j}}(k-1)+\widehat{\psi_{\ast j}}(k)}}{2} \bigl(\widehat{\psi_{\ast j}}(k-1)-\widehat{\psi_{\ast j}}(k)\bigr)+
 \int_{\r} \widehat{f_{\ast j}}'(\xi)\overline{\widehat{f_{\ast j}}(\xi)}\,d\xi
$$
$$
= \sum_{k\in\z} \int_{q_j^{-1}(k-1)}^{q_j^{-1} k}
\left(
\frac{\overline{\widehat{f_{\ast j}}\bigl(q_j^{-1}(k-1)\bigr)+\widehat{f_{\ast j}}\bigl(q_j^{-1}(k)\bigr)}}{2}\frac{\widehat{f_{\ast j}}\bigl(q_j^{-1}(k-1)\bigr)-\widehat{f_{\ast j}}\bigl(q_j^{-1}(k)\bigr)}{q_j^{-1}}+\widehat{f_{\ast j}}'(\xi)\overline{\widehat{f_{\ast j}}(\xi)}
\right)\,d\xi 
$$
The function $\widehat{f_{\ast j}}$ is linear on $[q_j^{-1}(k-1),\,q_j^{-1} k]$, therefore,  
$ \widehat{f_{\ast j}}'(\xi) \equiv -\left(\widehat{f_{\ast j}}\bigl(q_j^{-1}(k-1)\bigr)-\widehat{f_{\ast j}}\bigl(q_j^{-1}(k)\bigr)\right)q_j$ 
 on $[q_j^{-1}(k-1),\,q_j^{-1} k],$ 
denoting $c_{j,k}:=-\left(\widehat{f_{\ast j}}\bigl(q_j^{-1}(k-1)\bigr)-\widehat{f_{\ast j}}\bigl(q_j^{-1}(k)\bigr)\right)q_j$ and continuing calculations, we obtain
$$
\sum_{k\in\z} c_{j,k} \int_{q_j^{-1}(k-1)}^{q_j^{-1} k}\left( - \frac{\overline{\widehat{f_{\ast j}}\bigl(q_j^{-1}(k-1)\bigr)+\widehat{f_{\ast j}}\bigl(q_j^{-1}(k)\bigr)}}{2} +\overline{\widehat{f_{\ast j}}(\xi)}
\right)\,d\xi = 0,
$$
where the last equality is again due to the linearity of $\widehat{f_{\ast j}}$ on $[q_j^{-1}(k-1),\,q_j^{-1} k].$
So items 1.---5. are estimated.

Now, we write the squared Breitenberger UC in the form 
$$
UC_B^2(\psi_{\ast j}) = q_j^{2} {\rm var_A } (\psi_{\ast j}) 
q_j^{-2} {\rm var_F } (\psi_{\ast j}),
$$  
where 
$$
q_j^{-2} {\rm var_F } (\psi_{\ast j}) = q_j^{-2} \left(\frac{\|\psi_{\ast j}'\|^2}{\|\psi_{\ast j}\|^2}+\frac{(\psi_{\ast j}',\, \psi_{\ast j})^2}{\|\psi_{\ast j}\|^4}\right)
$$
and  (see also (\ref{defAB}), (\ref{AB}))
$$
q_j^{2} {\rm var_A } (\psi_{\ast j}) = q_j^{2}\left(\frac{\|\psi_{\ast j}\|^4}{|\tau(\psi_{\ast j})|^2}-1\right)=
q_j^{2} \frac{\|\psi_{\ast j}\|^4-|\tau(\psi_{\ast j})|^2}{|\tau(\psi_{\ast j})|^2}
=
q_j^{2} \frac{2 A(\psi_{\ast j}) \|\psi_{\ast j}\|^2- A^2(\psi_{\ast j})+B^2(\psi_{\ast j})}{\bigl(\|\psi_{\ast j}\|^2-A(\psi_{\ast j})\bigr)^2-B^2(\psi_{\ast j})}.
$$
Using (\ref{cond3}), (\ref{cond4}), and $\|\psi_{\ast j}\|^2 - \|\psi_{j}\|^2= \|\psi_{\ast j}\|^2 - 1 \to 0$ as $j \to \infty$, we see that 
$q_j^{-2} {\rm var_F } (\psi_{\ast j})$ is bounded as $j\to \infty.$ And since, by (\ref{UCBpsi0j}), there exists finite $\lim_{j\to \infty}UC_B(\psi_{\ast j}),$ it follows that there exists an absolute constant $C_0>0$ such that $q_j^{2} {\rm var_A } (\psi_{\ast j})> C_0$ as $j\to \infty$. 
Similarly, (\ref{cond5}), (\ref{cond6}), and $\|\psi_{\ast j}\|^2 - \|\psi_{j}\|^2= \|\psi_{\ast j}\|^2 - 1 \to 0$ as $j \to \infty$, yield the boundedness of $q_j^{2} {\rm var_A } (\psi_{\ast j})$  as $j\to \infty.$ Therefore,  there exists an absolute constant $C_0>0$ such that $q_j^{-2} {\rm var_F } (\psi_{\ast j})> C_0$ as $j\to \infty$. 

Hence the inequality $q_j^{2} {\rm var_A } (\psi_{\ast j})> C_0>0$ as $j\to \infty$ and estimations of items 1.,4.,5. enable to write
$$
\frac{\Delta^2(f_{\ast j})}{q_j^{2} {\rm var_A }(\psi_{\ast j})}=\frac{\|f_{\ast j}\|^2 \|\cdot f_{\ast j}\|^2-(\cdot f_{\ast j},\,f_{\ast j})^2}{\|f_{\ast j}\|^4}
\frac{(\|\psi_{\ast j}\|^2-A(\psi_{\ast j}))^2-B^2(\psi_{\ast j})}{q_j^{2}\left(2 A(\psi_{\ast j}) \|\psi_{\ast j}\|^2-A^2(\psi_{\ast j})+ B^2(\psi_{\ast j})\right)}
$$ 
$$
=
\frac{\left(\|\psi_{\ast j}\|^2+o(1)\right) 2 \cdot q_j^{2} A(\psi_{\ast j}) + q_j^{2} B^2(\psi_{\ast j})}{q_j^{2}\left(2 A(\psi_{\ast j}) \|\psi_{\ast j}\|^2-A^2(\psi_{\ast j})+ B^2(\psi_{\ast j})\right)}
\frac{\left((\|\psi_{\ast j}\|^2-A(\psi_{\ast j}))^2-B^2(\psi_{\ast j})\right)}{(\|\psi_{\ast j}\|^2+o(1))^2}
$$
$$
=\left(1+\frac{o(1) 2 \cdot q_j^{2} A(\psi_{\ast j}) + q_j^{2} A^2(\psi_{\ast j})}{q_j^{2}\left(2 A(\psi_{\ast j}) \|\psi_{\ast j}\|^2-A^2(\psi_{\ast j})\right)}\right)
\frac{\left((\|\psi_{\ast j}\|^2-A(\psi_{\ast j}))^2-B^2(\psi_{\ast j})\right)}{(\|\psi_{\ast j}\|^2+o(1))^2} \rightarrow 1 \mbox{ as } j\to \infty
$$
Similarly,  by the inequality $q_j^{-2} {\rm var_F } (\psi_{\ast j})> C_0 > 0$ as $j\to \infty$ and estimations of items 1.,2.,3., we conclude that
$$
\frac{\Delta^2(\widehat{f_{\ast j}})}{q_j^{-2} {\rm var_F }(\psi_{\ast j})}=
\frac{\|\psi_{\ast j}\|^4 \left(\|i f'_{\ast j}\|^2 \|f_{\ast j}\|^2 + (f'_{\ast j},\,f_{\ast j})^2\right)}{\|f_{\ast j}\|^4\  q_j^{-2}\left(\|\psi'_{\ast j}\|^2 \|\psi_{\ast j}\|^2 + (\psi'_{\ast j},\,\psi_{\ast j})^2\right)} 
$$
$$
=\frac{\|\psi_{\ast j}\|^4 \left((q_j^{-2}\|\psi'_{\ast j}\|^2 + o(1))(\|\psi_{\ast j}\|^2+o(1))+(q_j^{-1}(\psi'_{\ast j},\,\psi_{\ast j})+o(1))^2\right)}{(\|\psi_{\ast j}\|^2+o(1))^2 \ q_j^{-2}\left(\|\psi'_{\ast j}\|^2 \|\psi_{\ast j}\|^2 + (\psi'_{\ast j},\,\psi_{\ast j})^2\right)} \rightarrow 1 
\mbox{ as } j \to \infty.
$$
Finally, (\ref{UCHUCB}) follows from the last two limits and existence of finite $\lim_{j\to \infty}UC_B(\psi_{\ast j}).$

\textbf{Step 4.}
Let us consider auxiliary functions $f_j$, $j \in \n$ defined by
$$
\widehat{f_{j}} := \widehat{f_{\ast j}}(\cdot + c(\widehat{f_{\ast j}})),
$$ 
where $c(\widehat{f})$ is the frequency centre of the function $f$ (see Definition \ref{Huc}). Then it is well-known (see \cite[Exercise 1.5.1]{NPS} ) that $c(\widehat{f_j})=0$ and 
\be
\label{UCH0}
UC_H(f_{\ast j})=UC_H(f_j),\ \ \ \ \ \ \ \ 
 j \in \n.
\ee 
Let us check that $|c(\widehat{f_{\ast j}})|\leq q_j^{-1} M(C)$ as $j>j_0$ for some $j_0\in\n.$ 
Indeed,
by definition of a frequency centre, estimations (\ref{estnorm}), (\ref{estcentre}) yields 
$$
|c(\widehat{f_{\ast j}})| = \frac{\left(\cdot \widehat{f_{\ast j}},\,\widehat{f_{\ast j}} \right)}{\|\widehat{f_{\ast j}}\|^2}  
=
\frac{\left|\left(f'_{\ast j},\,f_{\ast j}\right)\right|}{\|f_{\ast j}\|^2}
\leq 
 q_j^{-1} \frac{\left|(\psi'_{\ast j},\,\psi_{\ast j})\right|+C \sqrt{2C}/3 + 1/6}{\|\psi_{\ast j}\|^2-C q_j^{-2}/3}
$$ 
Since  $\|\psi_{\ast j}\|^2-C q_j^{-2}/3 = 1 -C q_j^{-2}/3  +o(1) \geq 1/2$ as $j>j_0$ for some $j\in\n$, it follows from (\ref{cond3}) and definition of $M(C)$ (\ref{cond1}) that
$$
|c(\widehat{f_{\ast j}})| \leq 2\cdot q_j^{-1}\left(\left|(\psi'_{\ast j},\,\psi_{\ast j})\right|+C \sqrt{2C}/3 + 1/6\right) 
< 2\cdot q_j^{-1}\left(C+C \sqrt{2C}/3 + 1/6\right)
=q_j^{-1} M(C).
$$
Finally, using conditions (\ref{condpsi0}) and (\ref{conf1}) we conclude that
$\widehat{f_{\ast j}}(\xi) \equiv 0$ for $\xi \in [ - q_j^{-1}M(C),\,q_j^{-1}M(C)].$
Therefore, the inequality $|c(\widehat{f_{\ast j}})|\leq q_j^{-1} M(C)$ as $j>j_0$ for some $j_0\in\n$ provides  $\widehat{f_j}(0)=\widehat{f_{\ast j}}(c(\widehat{f_{\ast j}}))=0$ 
as $j>j_0$ for some $j_0\in\n$. 

\textbf{Step 5.}
Since $\widehat{f_j}(0)=0$ as $j>j_0$ for some $j_0\in\n$ and $c(\widehat{f_{j}})=0$ it follows from Theorem \ref{Battle}
that $UC_H(f_j)\geq 3/2$ as $j>j_0$ for some $j_0\in\n$. It remains to note that using (\ref{UCBpsi0j}), (\ref{UCHUCB}), and (\ref{UCH0}) we obtain 
$$
\lim_{j\to \infty} UC_B(\psi_j)= \lim_{j\to \infty} UC_B(\psi_{\ast j})=\lim_{j\to \infty} UC_H(f_{\ast j})=\lim_{j\to \infty} UC_H(f_j)\geq 3/2.
$$
Theorem \ref{main} is proved. \hfill $\Diamond$

\vspace{0.5cm}

Analyzing the conditions of Theorem \ref{main} one could ask a natural question about classes of periodic sequences that satisfy these conditions. Are these conditions restrictive or mild?  Let us make some illuminating remarks. 

\vspace{0.5cm}

\textbf{Remark 1: Wavelet sequence generated by periodization}

\noindent
Let $\psi^0 \in L_2(\mathbb{R})$ be a wavelet function on the real line. It is natural to require that $\cdot \psi^0,\, i (\psi^0)' \in L_2(\r)$. Otherwise $UC_H(\psi^0)$ will be infinite.   
Put 
$$
\psi^p_{j,k}(x):=2^{j/2}\sum_{n \in \mathbb{Z}} \psi^0(2^j(x+2 \pi n)+k).
$$
{\sloppy
The sequence $(\psi^p_{j,k})_{j,k},$ $j \in \z_+,$ $k=0,\dots,2^j-1$ is said to be \texttt{a periodic wavelet set generated by periodization}. Set 
$\psi^p_{j} := \psi^p_{j,0},$ $j = 0,1,\dots.$ Put $q_j=2^j.$ We claim that for a wavelet sequence $(\psi^p_{j})_{j\in\z_+}$ conditions (\ref{cond4})-(\ref{cond6}) are fulfilled, and $\lim_{j\to\infty}UC_B(\psi^p_j)$ is finite. If additionally $\widehat{\psi^0}(\xi) = o (\sqrt{\xi})$ as $\xi 
\to 0$, and $\cdot  (\psi^0)' \in L_2(\r),$ then condition (\ref{cond1}), and (\ref{cond2}) are also fulfilled respectively. Indeed, in \cite{prqurase03}, it is proved that the quantities 
$$
\|\psi^p_j\|_{L_{2,2\pi}}, \ \ 
2^{-2j}\|(\psi^p_j)'\|^2_{L_{2,2\pi}},\ \ 
2^{-j}((\psi^p_j)',\,\psi^p_j)_{L_{2,2\pi}},
\ \ 
2\cdot 2^{2j} A(\psi^p_j), 
\mbox{ and }\  2^j B(\psi^p_j)
$$ 
tend to 
$$
\|\psi^0\|_{L_2(\r)},\ \ 
\|(\psi^0)'\|^2_{L_2(\r)},\ \ 
((\psi^0)',\,\psi^0)_{L_2(\r)},
\ \ 
\|\cdot \psi^0\|^2_{L_2(\r)}, \mbox{ and }  
\ i(\cdot \psi^0,\,\psi^0)_{L_2(\r)}
$$
 respectively. Therefore, conditions (\ref{cond4}) - (\ref{cond6}) hold true. Then in \cite{prqurase03} it is deduced that 
  $\lim_{j\to\infty}UC_B(\psi^p_j) = UC_H(\psi^0),$ so   $\lim_{j\to\infty}UC_B(\psi^p_j)$ is finite. Since $\widehat{\psi^p_j}(k) = 2^{-j/2}\widehat{\psi^0}(2^{-j}k)$ and  $\widehat{\psi^0}(\xi) = o (\sqrt{\xi})$ as $\xi 
\to 0$, it follows that (\ref{cond1}) is satisfied for any $k \in \n.$
The condition 
$\widehat{\psi^0}(\xi) = o (\sqrt{\xi})$ as $\xi 
\to 0$  is not restrictive.
If this condition 
is not fulfilled, then $UC_H(\psi_0)=\infty.$  Indeed, suppose $UC_H(\psi_0)$ is finite. 
Then $(1+|\cdot|)\psi,\, (1+|\cdot|)\widehat{\psi}  \in L_2(\r),$ and $\psi^0$ is absolutely continuous (see \cite[Theorem 1.5.2]{NPS}). Therefore, $\psi^0(x) = O (x^{-3/2-\varepsilon}),$ $\varepsilon > 0$ 
as $x\to\infty.$ 
Then $|\widehat{\psi^0}(\xi+h) - \widehat{\psi^0}(\xi)| = O(h^{1/2+\varepsilon/2})$ as $\xi\in\r.$   And the condition 
$\widehat{\psi^0}(\xi) = o (\sqrt{\xi})$ as $\xi 
\to 0$ follows from the fact that $\widehat{\psi}(0)=0.$    
Finally, $\lim_{j\to \infty} 2 A((\psi^p_j)') = \|\cdot (\psi^0)'\|_{L_2(\r)}$ and 
$\cdot (\psi^0)' \in L_2(\r)$ yields (\ref{cond2}).

}

It is clear that combining the result $\lim_{j\to\infty}UC_B(\psi^p_j) = UC_H(\psi^0)$ of \cite{prqurase03} and Theorem \ref{Battle} we immediately get the inequality $\lim_{j \to \infty} UC_B(\psi^p_j)\geq 3/2$, and to do so  we do not need condition (\ref{cond2}). However the above inequality is proved here  under the condition
$((\psi^p_j)',\,\psi^p_j)_{L_{2,2\pi}} =0,$ while   only the  mild restriction
 $((\psi^p_j)',\,\psi^p_j)_{L_{2,2\pi}}\leq C \|\psi^p_j\|^2_{L_{2,2\pi}}$ 
 is required in Theorem \ref{main}, cf. (\ref{cond3}).     
 
\vspace{0.5cm}

\textbf{Remark 2: Wavelet sequence generated by UEP.}

\noindent
Let $(\psi_j)_{j\in\z_+}$ be a wavelet sequence satisfying Theorem  \ref{UEP} and $c_1 \leq \|\psi_j\|\leq c_2$, where $c_1>0,$ $c_2$ are absolute constants. Put $q_j=2^{j/2}.$ Then  condition (\ref{cond1})  is fulfilled for $k \in \z$. In fact, 
by (\ref{con3}) and (\ref{con4}), we conclude that
$
\widehat{\psi_j}(k)= e^{2\pi i 2^{-j-1} k} \mu_{k+2^j}^{j+1} \widehat{\varphi_{j+1}}(k)
$ 
and it follows from (\ref{con4}) 
$
|\mu_{k+2^j}^{j+1}|\leq \sqrt{2}. 
$
Therefore, 
$
|\widehat{\psi_j}(k)|\leq \sqrt{2} |\widehat{\varphi_{j+1}}(k)|.
$
Using 
(\ref{con1}), we get (\ref{cond1}) for all $k \in \n$.
Next, if $\psi_j$ is a trigonometric polynomial of degree less than  $C_1 2^{j/2}$, where $C_1$ is an absolute constant,  then condition (\ref{cond4}) is also fulfilled. Indeed, it follows from the Bernstein inequality that   
$
2^{-j}\|\psi_j'\|^2 / \|\psi_j\|^2 \leq C^2_0.
$

\vspace{0.5cm}

\textbf{Remark 3: Wavelet sequence constructed in \cite{LebPres14}.}

\noindent
Let $(\psi^a_j)_{j\in\n},$ $a>1,$ be a family of periodic wavelet sequences constructed in \cite{LebPres14}. Then all conditions (\ref{cond1})--(\ref{cond6}) are hold true for $q_j = j^{1/2}.$ Indeed, it is straightforward to see that 
$\widehat{\psi^a_j}(k) \asymp j^{-1} 2^{-j/2}$ for a fixed $k \in \z,$
$\|\psi^a_j\|^2 \asymp j^{-1/2} 2^{-j},$
$A((\psi^a_j)') \asymp j^{-1/2} 2^{-j}$,
$\|(\psi^a_j)'\|^2 \asymp j^{1/2} 2^{-j},$
$A(\psi^a_j) \asymp j^{-3/2} 2^{-j},$
$|B(\psi^a_j)| \asymp \sin(2 \pi 2^{-j-1}) j^{-1/2} 2^{-j}$
as $j \to \infty.$

\section{On a particular minimization problem for the Heisenberg UC}
One could make an analogy between Theorem \ref{Battle} and Theorem \ref{main}. For $f\in L_2(\r)$, exact conditions $c(\widehat{f})=0$ and $\int_{\r} f = \widehat{f}(0) = 0$ of Theorem  \ref{Battle} correspond to mild limit conditions (\ref{cond3}) and (\ref{cond1}) respectively. So it would be expectable to get a generalization of Theorem \ref{Battle} of the following  form: if $\int_{\r} f = \e$ and $c(\widehat{f})=0$, then $UC_H(f)\geq \alpha(\e),$ where $\lim_{\e\to 0}\alpha(\e)=3/2.$   Unfortunately, it is impossible to generalize the proof of Theorem \ref{Battle} to the above case. It turns out (Theorem \ref{nominfunc}) that there is no function satisfying the following minimization problem 
$$
UC_H(f) \to \min;\  f \in L_2(\r),\  
\int_{\r} f = \varepsilon, \ \e >0,\   
c(\widehat{f})=0.
$$ 
In the case of Theorem \ref{Battle} ($\e=0$) such a function do exist.
 To suggest an alternative proof is an open question.  
 \begin{teo}
\label{nominfunc}
Let $f$ be a function such that $\ \cdot f, \ i f' \in L_2(\r),$ $f_0:=a^{1/2}f(a\cdot),$ where $a=(\|\cdot f\| / \|i f'\|)^{1/2}$, and $c(\widehat{f_0})=0$. Fix $\e\in\cn$ and suppose  $\int_{\r}f_0(x)\,dx=\e.$ Then if $\e\neq 0,$ $\e\neq 2 \pi^{3/4}\frac{\sqrt{(2k)!}}{2^k k!},$ $k\in\n,$ 
 then  
there is no function satisfying the following minimization problem
\be
\label{minprob}
\left\{
\begin{array}{l}
UC_H(f_0) \to \min,\\
\int_{\r} f_0 = \varepsilon, \ \ \varepsilon >0, \ \ \ c(\widehat{f_0})=0.
\end{array}
\right.
\ee
If 
$\e = 2 \pi^{3/4}\frac{\sqrt{(2k)!}}{2^k k!},$ $k\in\n$,
then 
the Hermite function 
\be
\label{Hermite}
\phi_n(x)=\left(\frac{2^n n!}{2\sqrt{\pi}}\right)^{-1/2}(-1)^n e^{x^2/2} \frac{d^n}{dx^n}(e^{-x^2}), 
\ee
$n=2k,$ minimizes the problem (\ref{minprob}) and $UC_H(\phi_{2k})=(4k+1)/2$.
 \end{teo}

\textbf{Proof.}
The case $\e=0$ is considered in \cite[p.137-138]{Bat} (see Theorem \ref{Battle}). We exploit the ``variational'' idea  from the proof of this theorem for an arbitrary, sufficiently small parameter  $\e$.    
Suppose $f_1(x):=f_0(x+x_{0f}),$ then the time and frequency centres  of the  function $f_1$ are equal to zero, and $\Delta_{f_0}=\Delta_{f_1},$
$\Delta_{\widehat{f_0}}=\Delta_{\widehat{f_1}},$ $\int_{\r}f_1 = \int_{\r} f_0.$ It is also clear that $UC_H$ is a homogeneous functional that is $UC_H(c f)=UC_H(f)$ for $c\in \cn \setminus \{0\}$. So without loss of generality we assume that the functions $f_0$ and $\widehat{f_0}$ have zero centres and $\|f_0\|=1.$ Then $UC_H$ takes the form $UC_H(f_0)=\|\cdot f_0\| \|i f'_0\|.$ Let us minimize the functional  $\|\cdot f_0\|^2 \|i f'_0\|^2$ over functions $f_0\in L_2(\r)$ satisfying the constraints $\|f_0\|=1$ and 
$\int_{\r}f_0=\e.$ 
For the Lagrange function  
$$
\|\cdot f_0\|^2 \|i f'_0\|^2 + \lambda(\|f_0\|^2-1)+\kappa \left(\int_{\r}f_0 - \e\right)
$$
we get the Euler-Lagrange equation 
\begin{equation}
\label{dif_eq}
\|i f'_0\|^2 x^2 f_0(x) - \|\cdot f_0\|^2 f''_0(x) + \lambda f_0(x)+ \frac12 \kappa = 0.
\end{equation}
By definition of $f_0$ we get  
$\|i f'_0\| = a\|i f'\| = a^{-1}\|\cdot f\| = \|\cdot f_0\|=:\sqrt{\alpha/2}$. 
Then the equation is rewritten as
$$
\frac12 \alpha (x^2 f_0(x)-f''_0(x)) +\lambda f_0(x) +\frac12 \kappa = 0,
$$
where $\frac12 \alpha$ is a value of the desired functional.

If $\kappa=0$ then the solutions are eigenfunctions of the Hamilton operator $H f=\frac12 (\cdot^2 f-f'')$ that is the Hermite functions (\ref{Hermite})
normalized by the condition $\|\phi_n\|^2=(2\pi)^{-1}\int |\phi_n|^2 = 1.$  It can be proved by fairly straightforward calculation that 
$$
\int_{\r}\phi_n=\left\{
\begin{array}{ll}
0, & n=2k+1, \\
2 \pi^{3/4}\frac{\sqrt{(2k)!}}{2^k k!}, & n=2k.
\end{array}
\right.
$$
and $a=(\|\cdot \phi_n\| / \|i \phi'_n\|)^{1/2} = 1$.
If $\e\neq\int_{\r}\phi_n,$ then there are no solutions of minimization problem satisfying the constraints. Otherwise, the Hermite functions (\ref{Hermite}) minimizes the problem (\ref{minprob}) and $UC_H(\phi_n)= \|\cdot \phi_n\| \|i \phi_n'\| =(2n+1)/2.$ 

If $\kappa\neq 0,$ then we use the completeness of the set of Hermite functions, expand a function $f_0$ into a series $\sum a_n\phi_n,$ and substitute it in the equation. As a result we get 
$$
\sum_{k=0}^{\infty} a_k (\alpha(k+1/2)+\lambda)\phi_k(x)+ \kappa/2 = 0.
$$   
Orthonormality of the set $\{\phi_n\}_{n\in\n}$ allows to find the coefficients 
$$
a_k = - \frac12 \kappa \frac{1}{\alpha(k+1/2)+\lambda} \int_{\r} \phi_k.
$$ 
Multiplying  equation (\ref{dif_eq}) by $f_0$ and integrating over $\r$ we obtain a relationship between parameters
$\alpha^2/2+\lambda + \kappa \e/2 = 0. $
Therefore, the solution of (\ref{dif_eq}) takes the form
$$
f_0(x) = -\pi^{3/4} \frac{\kappa}{\alpha}\sum_{n=0}^{\infty} 
\frac{\sqrt{(2n)!}}{2^n n!}\frac{\phi_{2n}(x)}{2n+1/2-(\alpha+\kappa \e/\alpha)/2}.
$$
The constraints give us the following equations 
\begin{gather}
\e=\int_{\r}f_0 = - 2 \pi^{3/2} \frac{\kappa}{\alpha}\sum_{k=0}^{\infty}
\frac{(2n)!}{(2^n n!)^2}\frac{1}{2n+1/2-(\alpha+\kappa \e/\alpha)/2}
\nonumber
\\
1=\frac{1}{2\pi}\int_{\r}|f_0|^2 = \pi^{3/2}\frac{\kappa^2}{\alpha^2}\sum_{k=0}^{\infty}
\frac{(2n)!}{(2^n n!)^2}\frac{1}{(2n+1/2-(\alpha+\kappa \e/\alpha)/2)^2}
\nonumber
\end{gather}
Using \cite[Chap.12, Ex.8]{WW} we rewrite the equations
 \begin{gather}
\e= -  \pi^{3/2} \frac{\kappa}{\alpha}\Gamma\left(\frac12\right)\frac{\Gamma\left(\frac14-\frac14\left(\alpha+\frac{\kappa \e}{\alpha}\right)\right)}{\Gamma\left(\frac34-\frac14\left(\alpha+\frac{\kappa \e}{\alpha}\right)\right)}
\nonumber
\\
1= \frac{\pi^{1/2}}{2}\frac{\kappa^2}{\alpha^2}\Gamma\left(\frac12\right)
\frac{d}{d(\alpha+\kappa \e/\alpha)}\frac{\Gamma\left(\frac14-\frac14\left(\alpha+\frac{\kappa \e}{\alpha}\right)\right)}{\Gamma\left(\frac34-\frac14\left(\alpha+\frac{\kappa \e}{\alpha}\right)\right)}.
\nonumber
\end{gather}
We introduce notions $\beta:=\kappa/\alpha,$ $F(x):=\Gamma(1/2)
\frac{\Gamma(1/4-x/4)}{\Gamma(3/4-x/4)}.$ Then the last equations take the form
\begin{gather}
\nonumber
\left\{
\begin{array}{l}
-\beta F(\alpha+\beta \e) = \pi^{-3/2} \e,
 \\
\beta^2 F'(\alpha+\beta \e) = 2 \pi^{-1/2}.
\end{array}
\right.
\end{gather}
We claim that this system has no solutions if $\varepsilon\neq 0$. Indeed, if $\e = 0,$ then $\alpha = 3,$ $\beta^2 = 2 \pi^{-1/2} (F'(3))^{-1}.$ (This is the case of Theorem \ref{Battle}.) 
The function $F$ is continuously differentiable and increasing in the neighborhood of $3$ and $F(3)=0.$
However it follows from the equation 
$
-\beta F(\alpha+\beta \e) = \pi^{-3/2} \e
$
that $F(3+0)<0$ and $F(3-0)>0$. So $F$ can not be an increasing function.  This contradiction means that in the case $\kappa \neq 0,$ there is no solution of minimization problem (\ref{minprob}). Theorem \ref{nominfunc} is proved.
\hfill $\Diamond$

\section*{Acknowledgments}
The work is supported by the RFBR, grant \#15-01-05796, and 
 by Saint Petersburg State University, grant  \#9.38.198.2015.

\end{document}